\documentclass[12pt]{article}
\usepackage{latexsym}
\usepackage{amsmath, amsthm, amsfonts}
\usepackage{graphics}
\usepackage{graphicx}
\usepackage{color}
\usepackage{float}

\usepackage{array}
\usepackage{hyperref}

\newcommand{\abs}[1]{\left\lvert#1\right\rvert}
\newtheorem{theorem}{Theorem}[section]

\newtheorem{conjecture}{Conjecture}[section]

\def\smallskip{\addvspace{\smallskipamount}}
\def\medskip{\addvspace{\medskipamount}}
\def\bigskip{\addvspace{\bigskipamount}}
\def\makefootline{\baselineskip=24pt \line{\the\footline}}
\def\pagecontents{\ifvoid\topins\else\unvbox\topins\fi
   \dimen@=\dp255 \unvbox255
   \ifvoid\footins\else
      \vskip\skip\footins \footnoterule \unvbox\footins\fi
     \ifr@ggedbottom \kern-\dimen@ \vfil \fi}
\def\footnoterule{\kern-3pt\hrule width 2truein \kern 2.6pt}

\begin{document}

\title{{Dynamics of the complex rational delay recursive sequence $\displaystyle{z_{n+1}=\frac{\alpha+\beta z_{n-k}}{\gamma - z_{n}}}$}}

\small \author{Sk. Sarif Hassan\\
  \small {Department of Mathematics}\\
  \small {University of Petroleum and Energy Studies}\\
  \small {Bidholi, Dehradun, India}\\
  \small Email: {\texttt{\textcolor[rgb]{0.00,0.07,1.00}{s.hassan@ddn.upes.ac.in}}}\\
}

\maketitle
\begin{abstract}
\noindent Dynamics of the delay rational difference equation $\displaystyle{z_{n+1}=\frac{\alpha+\beta z_{n-k}}{\gamma - z_{n}}}$ with complex parameters $\alpha$, $\beta$, $\gamma$ and arbitrary complex initial conditions is investigated. Existence of prime period two solutions and higher order periods are ensured in the complex parameters unlike in the case of real parameters of the same rational difference equation. In addition, a new dynamical behavior, chaotic solutions of the difference equation are ensured computationally.
\end{abstract}

\begin{flushleft}\footnotesize
{Keywords: Rational difference equation, Local asymptotic stability, Chaotic trajectory and Periodicity. \\

}
\end{flushleft}

\section{Introduction and Background}

Consider the rational difference equation,

\begin{equation}
\displaystyle{z_{n+1}=\frac{\alpha+\beta z_{n-k}}{\gamma - z_{n}}},  n=0,1,2,\ldots
\label{equation:total-equationA}
\end{equation}%
where the parameters $\alpha$, $\beta$, $\gamma$ and the initial conditions $z_{-k},z_{-k+1},\dots z_{-1}, z_{0}$ are arbitrary complex numbers.

\addvspace{\bigskipamount}
\noindent
This (k+1)th order rational difference equation Eq.(\ref{equation:total-equationA}) and closely similar equations are studied when the parameters $\alpha \geq 0$, $\beta>0$ and $\gamma>0$ and the initial conditions are real numbers while $k$ is any natural number \cite{ZM}, \cite{XW} and \cite{CXW}. They have shown that the positive fixed point of Eq.(\ref{equation:total-equationA}) is a global attractor with a basin that depends on certain conditions posed on the coefficients $\alpha$, $\beta$ and $\gamma$ \cite{ZM}.

\noindent
In this present article, an attempt has been made to study the dynamics of the rational difference equation (\ref{equation:total-equationA}) under the assumption that the parameters and the initial conditions are arbitrary complex numbers. Similar works can be found in \cite{S-E}, \cite{S1} and \cite{S2}. Here some of the important basics are being discussed in the following \cite{S-H}, \cite{K-L} and \cite{Ko-L}:\\

\noindent

\textbf{Definition 1}: A difference equation of order $(k+1)$ is of the form

\begin{equation}
\displaystyle{z_{n+1}=f(z_n,z_{n-1},\dots,z_{n-k})},  n=0,1,2,\ldots
\label{equation:total-equationB}
\end{equation}%

\noindent
where $f$ is a continuous function which maps a subset $\mathbb{D}^{k+1}$ into $\mathbb{D}$ and $\mathbb{D}\subset\mathbb{C}$. A fixed point $\bar{z}$ of the difference equation Eq. (\ref{equation:total-equationB}) is a point that satisfy the condition $\bar{z}=f(\bar{z},\bar{z},\dots,\bar{z})$. \\

\textbf{Definition 2}: Let $\bar{z}$ be a fixed point of the Eq.(\ref{equation:total-equationB}), then $\bar{z}$ is \emph{locally asymptotically stable} if for every $\epsilon>0$, there exist a $\delta(\epsilon)>0$ such that, if $z_{-k},\dots,z_{-1},z_0 \in$ $D$ with $\abs{z_{-k}-\bar{z}}+\dots+\abs{z_{-1}-\bar{z}}+\abs{z_{0}-\bar{z}}<\delta(\epsilon)$, then $\abs{z_{n}-\bar{z}}<\epsilon$ for all $n \geq -k$.\\

\textbf{Definition 3}: A sequence ${z_n}_{n=−k}^{\infty}$ is said to be \emph{periodic with period $p$} if $z_{n+p} = z_n$ for all $n\geq-k$. A sequence ${z_n}_{n=−k}^{\infty}$ is said to be \emph{periodic} with prime period \emph{p} if $p$ is the smallest positive integer having this property.\\

\textbf{Definition 4}: An open ball $B(a,r) \in \mathbb{C}$ is called an \emph{invariant open ball} of Eq.(\ref{equation:total-equationB}) if $z_{-k},\dots,z_{-1},z_0 \in B(a,r)$ then $z_n \in B(a,r)$ for all $n >0$. That is every solution of Eq.(\ref{equation:total-equationB}) with initial conditions in $B(a,r)$ remains in $B(a,r)$.\\

\textbf{Definition 5}: The difference equation Eq.(\ref{equation:total-equationB}) is said to be \emph{permanent and bounded} if there exist positive real numbers $M$ and $N$ with $0<M\leq N<\infty$ such that for any initial conditions $z_{-k},\dots,z_{-1},z_0$ there exist a positive integer P which depends on the initial conditions such that $M\leq \abs{z_n}\leq N$ for all $n\geq P$.\\

\noindent
The linearized equation associated with Eq.(\ref{equation:total-equationB}) about the equilibrium point $\bar{z}$ is $$y_{n+1}=\sum_{i=0}^{k} \frac{\partial f(\bar{z},\bar{z},\dots,\bar{z})}{\partial u_i} y_{n-i}$$
\noindent
Its characteristic equation is $$\lambda^{k+1}=\sum_{i=0}^{k} \frac{\partial f(\bar{z},\bar{z},\dots,\bar{z})}{\partial u_i}\lambda^{n-i}$$
where $n=0,1,2,\dots$.

\begin{theorem}
Assume that $f$ is a $C^1$-function and let $\bar{z}$ a fixed point of Eq.(\ref{equation:total-equationB}). Then the following statements are true:
\begin{itemize}
  \item If all the roots of the characteristic equation lie in the open unit disk $\abs{\lambda}< 1$, then the fixed point $\bar{z}$ of Eq.(\ref{equation:total-equationB}) is locally asymptotically stable.
  \item If at least one root of the characteristic equation has the absolute value greater than one, then the fixed point $\bar{z}$ of Eq.(\ref{equation:total-equationB}) is unstable.
  \item If all the roots of the  have the characteristic equation absolute value greater than one, then the fixed point $\bar{z}$ of Eq.(\ref{equation:total-equationB}) is a source.
\end{itemize}

\end{theorem}

\begin{theorem}
  Assume that $p$, $q$ $\in \mathbb{C}$ and $k \in \mathbb{N}$. Then $\abs{p} + \abs{q} < 1$ is a
sufficient condition for asymptotically stability of the difference equation $$z_{n+1}-pz_n+qz_{n-k}=0, n=0,1,2,3,\dots$$
\end{theorem}

\noindent
Now we shall use these basic theorems to explore the local stability of the fixed points of the Eq.(\ref{equation:total-equationA}).

\section{Local Asymptotic Stability of the Fixed Points and Boundedness}
The fixed points of Eq.(\ref{equation:total-equationA}) are the solutions of the quadratic equation
\[
\bar{z}=\frac{\alpha+\beta \bar{z}}{\gamma-\bar{z}}
\]
Eq.(\ref{equation:total-equationA}) has the two fixed points $\bar{z}_1$ and $\bar{z}_2$ when $\abs{\alpha} \neq \abs{\frac{({\gamma-\beta})^2}{4}}$ \dots \\
$\frac{1}{2} \left(-\sqrt{-4 \alpha +\beta ^2-2 \beta  \gamma +\gamma ^2}-\beta +\gamma \right)$ and $\frac{1}{2} \left(\sqrt{-4 \alpha +\beta ^2-2 \beta  \gamma +\gamma ^2}-\beta +\gamma \right)$ respectively.\\

\noindent
It is noted that if $\alpha=\frac{({\gamma-\beta})^2}{4}$, then there is only one fixed point, $\frac{(\gamma-\beta)}{2}$. \\

\noindent
The linearized equation of the rational difference equation Eq.(\ref{equation:total-equationA}) with respect to the fixed points $\bar{z}$ is

\begin{equation}
\label{equation:linearized-equation}
\displaystyle{
z_{n+1} - \frac{\alpha +\beta  \bar{z}}{(\gamma -\bar{z})^2} z_{n} - \frac{\beta }{\gamma -\bar{z}} z_{n-k}=0,  n=0,1,\ldots
}
\end{equation}

\noindent
with associated characteristic equation

\begin{equation}
\lambda^{k+1} -\frac{\alpha +\beta  \bar{z}}{(\gamma -\bar{z})^2} \lambda^k -\frac{\beta }{\gamma -\bar{z}} = 0.\\
\end{equation}

\noindent
The linearized equation for the fixed point $\frac{(\gamma-\beta)}{2}$ is

\begin{equation}
\label{equation:linearized-equation}
\displaystyle{
z_{n+1} - \frac{\gamma-\beta}{\gamma+\beta} z_n - \frac{2 \beta }{\beta +\gamma } z_{n-k}=0,  n=0,1,\ldots
}
\end{equation}

\noindent
with associated characteristic equation

\begin{equation}
\lambda^{k+1} -\frac{\gamma-\beta}{\gamma+\beta} \lambda^k -\frac{2 \beta }{\beta +\gamma } = 0.
\end{equation}

\noindent
It is found that there does not exist any $\alpha$, $\beta$ and $\gamma$ in $\mathbb{C}$ such that $\alpha=\frac{({\gamma-\beta})^2}{4}$ in $\mathbb{D} \in \mathbb{C}$ for which the condition $$\abs{\frac{\gamma-\beta}{\gamma+\beta}}+\abs{\frac{2 \beta }{\beta +\gamma }}<1$$ Hence the fixed point $\frac{(\gamma-\beta)}{2}$ is not \emph{locally asymptotically stable}. \\

\noindent
The following result gives the local asymptotic stability of the fixed point \dots \\  $\frac{1}{2} \left(-\sqrt{-4 \alpha +\beta ^2-2 \beta  \gamma +\gamma ^2}-\beta +\gamma \right)$ and $\frac{1}{2} \left(\sqrt{-4 \alpha +\beta ^2-2 \beta  \gamma +\gamma ^2}-\beta +\gamma \right)$ of the Eq.(\ref{equation:total-equationA}).

\begin{theorem}
The fixed point $\frac{1}{2} \left(-\sqrt{-4 \alpha +\beta ^2-2 \beta  \gamma +\gamma ^2}-\beta +\gamma \right)$ of Eq.(\ref{equation:total-equationA}) is  locally asymptotically stable if $$\abs{\frac{2 \beta }{\sqrt{(\beta -\gamma )^2-4 \alpha }+\beta +\gamma }}+\abs{\frac{\gamma  \left(-\sqrt{(\beta -\gamma )^2-4 \alpha }-\beta +\gamma \right)-2 \alpha }{2 (\alpha +\beta  \gamma )}}<1$$
\end{theorem}

\begin{theorem}
The fixed point $\frac{1}{2} \left(\sqrt{-4 \alpha +\beta ^2-2 \beta  \gamma +\gamma ^2}-\beta +\gamma \right)$  of Eq.(\ref{equation:total-equationA}) is  locally asymptotically stable if $$ \abs{\frac{2 \beta }{-\sqrt{(\beta -\gamma )^2-4 \alpha }+\beta +\gamma }}+\abs{\frac{\gamma  \left(\sqrt{(\beta -\gamma )^2-4 \alpha }-\beta +\gamma \right)-2 \alpha }{2 (\alpha +\beta  \gamma )}}<1$$
\end{theorem}

\noindent
Proof of these two theorems are straightforward from the result stated in the \textbf{Theorem 1.2}. Here we go with few examples which illustrate the asymptotic behavior of these two fixed points.

\begin{table}[H]
\begin{tabular}{| m{1cm} || m{4cm} || m{9.5cm} |}
\hline
\centering  \textbf{S.n.} &
\begin{center}
\textbf{Parameters} $\alpha$, $\beta$, $\gamma$
\end{center}
 &
 \begin{center}
 \textbf{Remark}
 \end{center}
 \\
\hline
\centering  1 &
\begin{center}
$\alpha \to -107.7 - 374.6i, \beta \to 56 - 91.2i,\gamma \to 147.4 + 210i$
\end{center} & Here modulus of the zeros of the characteristic equation are both less than one and so the fixed point is \emph{locally asymptotically stable}. It is found that for all $k=1,2,3,\dots $, The trajectories are convergent to the fixed point $\frac{1}{2} \left(-\sqrt{-4 \alpha +\beta ^2-2 \beta  \gamma +\gamma ^2}-\beta +\gamma \right)$. \\
\hline
\centering  2 &
\begin{center}
$\alpha \to -297.4-75.3i, \beta \to 78-71.6i,\gamma \to -220.6 - 128.5i$
\end{center} & Here modulus of the zeros of the characteristic equation are less than one and so the fixed point $\frac{1}{2} \left(\sqrt{-4 \alpha +\beta ^2-2 \beta  \gamma +\gamma ^2}-\beta +\gamma \right)$ is \emph{locally asymptotically stable}. It is found that for all $k=1,2,3,\dots $, The trajectories are convergent to the fixed point.\\
\hline

\end{tabular}
\caption{Parameters $\alpha$, $\beta$ and $\gamma$ for which the fixed points $\frac{1}{2} \left(\pm \sqrt{-4 \alpha +\beta ^2-2 \beta  \gamma +\gamma ^2}-\beta +\gamma \right)$ of Eq.(\ref{equation:total-equationA}) are \emph{locally asymptotically stable} for arbitrary initial values.}
\label{Table:}
\end{table}

\begin{figure}[H]
      \centering

      \resizebox{16.5cm}{!}
      {
      \begin{tabular}{c c}
      \includegraphics [scale=9]{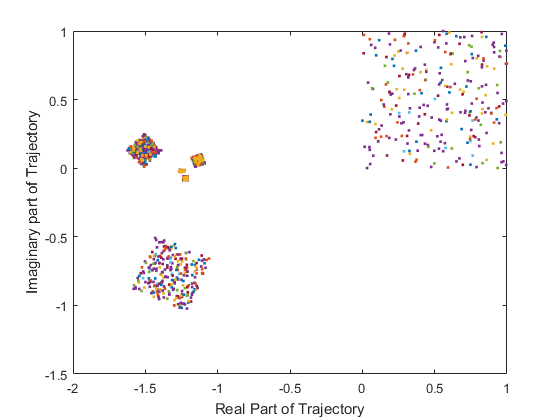}
      \includegraphics [scale=9]{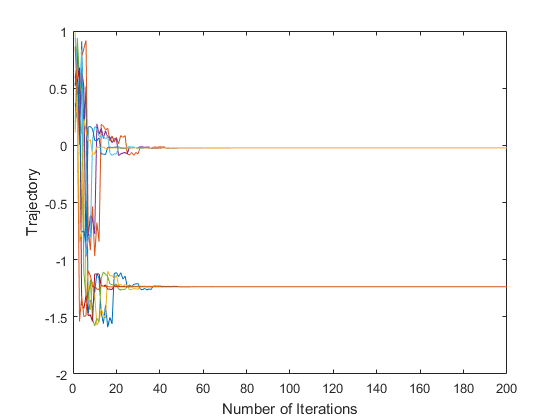}\\
      \includegraphics [scale=9]{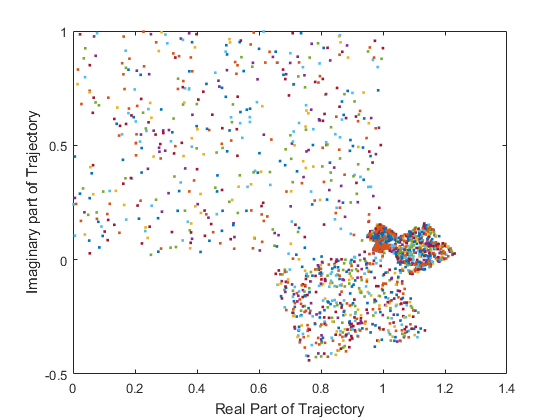}
      \includegraphics [scale=9]{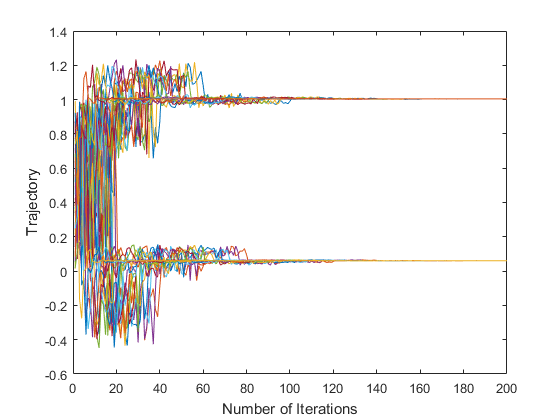}\\

\end{tabular}
      }
\caption{Local Asymptotical Stability of the fixed points. In the top of the figure, for different values of $k$, complex trajectories plot (left) and their corresponding time series plot with real and imaginary part (right) are given for the fixed point $\frac{1}{2} \left(-\sqrt{-4 \alpha +\beta ^2-2 \beta  \gamma +\gamma ^2}-\beta +\gamma \right)$. Same is given in the bottom for the other fixed point.}

      \end{figure}

\noindent
In the Table-1 as stated, for the parameters given in the serial number $1$ and $2$, the trajectories are convergent and converge to the fixed points $\frac{1}{2} \left(-\sqrt{-4 \alpha +\beta ^2-2 \beta  \gamma +\gamma ^2}-\beta +\gamma \right)$ and $\frac{1}{2} \left(\sqrt{-4 \alpha +\beta ^2-2 \beta  \gamma +\gamma ^2}-\beta +\gamma \right)$ of Eq.(\ref{equation:total-equationA}) respectively for all delay term $k$. \\ \\

\noindent
Further, a matrix of $20$ rows of parameters $\alpha$, $\beta$ and $\gamma$ for which the fixed point $\frac{1}{2} \left(-\sqrt{-4 \alpha +\beta ^2-2 \beta  \gamma +\gamma ^2}-\beta +\gamma \right)$ of Eq.(\ref{equation:total-equationA}) is  \emph{locally asymptotically stable}, is given here. In each row the first, second and third complex number referred as $\alpha$, $\beta$ and $\gamma$ respectively.

$$\left(
\begin{array}{ccc}
 2213.8\, -627.2 i & -1158.7-591.8 i & 4355.\, -33.7 i \\
 -1755.8-435.3 i & 1657.2\, -99.3 i & 3223.4\, -1018.8 i \\
 -880.8+1780.7 i & -1282.1-813.5 i & -555.-3870.8 i \\
 2860.1\, +1430 i & -1262.7-860 i & -643.8+3893.9 i \\
 -2501.1+1648.8 i & 1508.5\, +1652.8 i & 1773.1\, -4323.1 i \\
 -60.8-492 i & 7.9\, -90.7 i & 307.5\, +659 i \\
 2750.5\, -1175.1 i & -953.7+392.9 i & 3987.6\, +1450.2 i \\
 338.9\, +77.1 i & -203.4-84.8 i & 750.7\, +226.7 i \\
 -635.5+187.8 i & -224.9-840.4 i & 1347.8\, +1812.5 i \\
 18\, +69.1 i & 1.2\, +45.3 i & 276.5\, -26.2 i \\
 1055.8\, +3679.6 i & 501.5\, -369.8 i & 1400.3\, +1747.7 i \\
 2337.3\, +4233.9 i & -377.3-451.5 i & 120.7\, +1373.2 i \\
 3231.5\, +700.7 i & -968.6-1000.5 i & 3070.4\, -968.5 i \\
 2021.8\, +1066.1 i & -361.+739 i & 3616.2\, +1451.6 i \\
 -1058.1+1749.3 i & 357.7\, +4.8 i & 788.9\, -3083. i \\
 -890.5+1361.9 i & -1421.6-670.8 i & 116.4\, -3391.3 i \\
 -3635.7-4440.8 i & -1692.9+1093 i & 3719.2\, -3308.7 i \\
 -76.2+2700.8 i & -733.+1500 i & 3174.1\, +1390 i \\
 -1883.8-4906.7 i & -1490.2+742.2 i & 2899.5\, -4163.9 i \\
 301\, +1028.6 i & -112.9-1099.1 i & 2190.1\, -3783.9 i \\
\end{array}
\right)$$

\noindent
Similarly, the fixed point $\frac{1}{2} \left(\sqrt{-4 \alpha +\beta ^2-2 \beta  \gamma +\gamma ^2}-\beta +\gamma \right)$  of Eq.(\ref{equation:total-equationA}) is  \emph{locally asymptotically stable} for the set of parameters which is given as matrix as follows:

$$
\left(
\begin{array}{ccc}
 2686.1\, -1250.3 i & 629.7\, -417 i & -1594.+4310.2 i \\
 1419.2\, -2938.1 i & 588.9\, -568.9 i & -2005.2+421 i \\
 763.8\, +3108i & -646.-533.3 i & -1073.2-3104.1 i \\
 2816.4\, -1018.7 i & -503.9-401.1 i & -2364.7-2712. i \\
 -106.2-138.8 i & 231.4\, +106.1 i & -631.5+854.9 i \\
 2608.3\, -3767.4 i & -363.6-366.9 i & -1909.8+299.2 i \\
 3435.3\, +2539.7 i & 535.3\, +20.1 i & -618.1+2127.1 i \\
 -984.2+955.5 i & 666.8\, -431.4 i & 26.\, -2208.2 i \\
 -2113.7+3096.8 i & 991.3\, +434.3 i & -326.8-3814.8 i \\
 -669.8-794.3 i & -86.8-1545.9 i & -3677.5+1627.6 i \\
 2663.8\, -2211i & 1366.2\, +1397.7 i & -3862.-2238.9 i \\
 575\, +221i & 1515.2\, -8.2 i & 1087.6\, +4131.8 i \\
 992.7\, -772i & 180.9\, -1155.2 i & -2175.5-885.1 i \\
 3337.7\, +2555i & -1931.1+111.4 i & -2914.1-4788 i \\
 982.5\, +492.9 i & 807.6\, +666.8 i & -3412.6-1495.6 i \\
 1908.8\, +3639.9 i & -527.4+1553.7 i & -3264.5+645.1 i \\
 1393.1\, -853.1 i & 261.1\, +1276.7 i & -2089.8-4087.5 i \\
 -1931.3-1221.9 i & 340.2\, -527.6 i & -1151.6-2448i \\
 -1319.1+724.3 i & 394.1\, +798.2 i & -4187.5-992.1 i \\
 1582.8\, +350.7 i & 251.5\, +37.2 i & -2431.8+892.6 i \\
\end{array}
\right)$$

\noindent
The $100$ set of such parameters $\alpha$, $\beta$ and $\gamma$ are plotted in the complex plane in the following Fig. 2.

\begin{figure}[H]
      \centering

      \resizebox{8.5cm}{!}
      {
      \begin{tabular}{c}
      \includegraphics [scale=5]{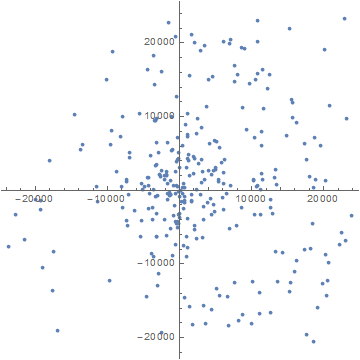}\\\\\\
      \includegraphics [scale=5]{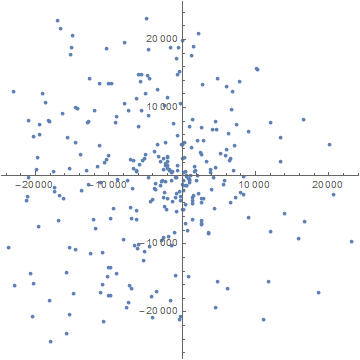}
\end{tabular}
      }
\caption{Parameters $\alpha$, $\beta$ and $\gamma$ plot for which the fixed points are locally asymptotically stable. Here the top figure stands for the fixed point $\frac{1}{2} \left(\sqrt{-4 \alpha +\beta ^2-2 \beta  \gamma +\gamma ^2}-\beta +\gamma \right)$ and bottom one stands for the fixed point $\frac{1}{2} \left(\sqrt{-4 \alpha +\beta ^2-2 \beta  \gamma +\gamma ^2}-\beta +\gamma \right)$.}

      \end{figure}

\subsection{Local Asymptotic Stability in the case $\alpha=0,  \in \mathbb{C}$}

When the parameter $\alpha =0$, then the the difference equation Eq.(\ref{equation:total-equationA}) would reduce to

\begin{equation}
\label{equation:linearized-equation11}
\displaystyle{
z_{n+1} =\frac{\beta z_{n-k}}{\gamma-z_n},  n=0,1,\ldots, k \in {1,2,\dots}
}
\end{equation}

\noindent
The fixed points $z_{\pm}$ of the equation Eq.(\ref{equation:linearized-equation11}) are $0$ and $\gamma-\beta$ respectively.

\noindent
In the similar fashion we did earlier, the fixed point $z_{+}=0$ of the difference equation Eq.(\ref{equation:linearized-equation11}) is \emph{locally asymptotically stable}, \emph{unstable} if $\abs{\frac{\beta}{\gamma}}<1$, $\abs{\frac{\beta}{\gamma}}>1$ respectively.\\

\noindent
Similarly, the fixed point $z_{-} =\gamma-\beta$ is \emph{locally asymptotically stable} if $\abs{\frac{\gamma}{\beta}-1}<0$ i.e. $\abs{\frac{\gamma}{\beta}}<1$. If $\abs{\frac{\gamma}{\beta}}>1$ then the fixed point is \emph{unstable}.

\subsection{Local Asymptotic Stability in the case $\alpha=\beta=\gamma$}

When all the three parameters are same $(\beta=\gamma=\alpha)$, then the difference equation Eq.(\ref{equation:total-equationA}) becomes

\begin{equation}
\displaystyle{z_{n+1}=\frac{\alpha+\alpha z_{n-k}}{\alpha - z_{n}}},  n=0,1,2,\ldots
\label{equation:total-equationD}
\end{equation}%

\noindent
The fixed points of the difference equation Eq.(\ref{equation:total-equationD}) are $\pm i\sqrt{\alpha}$. \\

\noindent
The fixed points $\pm i\sqrt{\alpha}$ of Eq.(\ref{equation:total-equationD}) are \emph{locally asymptotically stable} if $$\abs{\frac{2}{1+\mp\frac{i}{\sqrt{\alpha }}}}<1$$

\noindent
Here we shall look for the subset $\mathbb{S}$ of $\mathbb{C}$ of the parameter $\alpha$ for which the above stated condition does hold good. Before seeking the set of our interest, let us try to envisage the region and we found it to be the following in the Fig. 3 where the real part and imaginary part of the parameter $\alpha$ are assumed to varied over the region $(-1, 1)\times(-1, 1)$. \\

\noindent
Now we can see what we are to find, i.e. dependence of the boundary of this set, i.e. we should find a few functions yielding $y$ as a function of $x$ on the boundary. where $\alpha=x+iy$. We have the following $x$ and $y$ which satisfy the boundary $\abs{\frac{2}{1+\mp\frac{i}{\sqrt{\alpha }}}}=1$.

\begin{figure}[H]
      \centering

      \resizebox{11cm}{!}
      {
      \begin{tabular}{c}
      \includegraphics [scale=8]{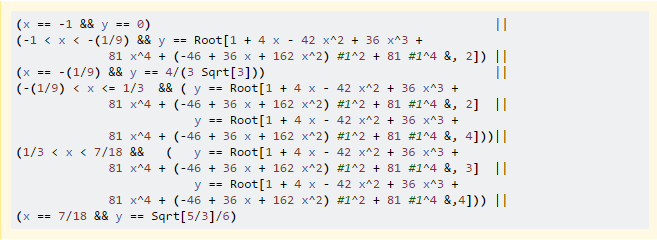}
      \end{tabular}
      }

      \end{figure}

\begin{figure}[H]
      \centering

      \resizebox{8cm}{!}
      {
      \begin{tabular}{c}
      \includegraphics [scale=5]{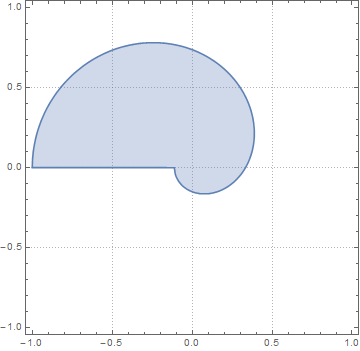}
      \end{tabular}
      }
\caption{Region of the parameter $\alpha$ while its real and imaginary part are varying along the (-1,1) and (-1,1) respectively.}

      \end{figure}

\noindent
To figure out why the solution looks slightly involved, we plot all the roots of the underlying polynomial then the curves will form more symmetric pattern:\\

\begin{figure}[H]
      \centering

      \resizebox{8cm}{!}
      {
      \begin{tabular}{c}
      \includegraphics [scale=5]{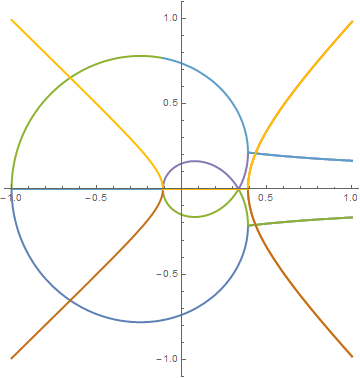}
      \end{tabular}
      }
\caption{Boundary of the roots of the underlying polynomials. Different colors denote four different polynomials.}

      \end{figure}

\noindent
Here we go with a set of examples of parameter $\alpha$ for which the fixed points $\pm i\sqrt{\alpha}$ of the difference equation Eq.(\ref{equation:total-equationD}) are \emph{locally asymptotically stable}. \\

\noindent
\emph{A=}$\{\{\alpha \to 0.294548\, +0.27933 i\},\{\alpha \to 0.303387 i\},\{\alpha \to -0.111111\},\{\alpha \to -0.111111+0.288786 i\},\{\alpha \to -0.12042 i\}\}$ \\

\noindent
\emph{B=}$\{\{\alpha \to -0.428571 i\},\{\alpha \to -0.526596-0.0285714 i\},\{\alpha \to -0.0725055\},\{\alpha \to 0.325133\, -0.0233645 i\},\{\alpha \to 0.28391\, +0.0116959 i\}\}$ \\ \\
\noindent
For these two sets \emph{A} and \emph{B} of parameters the fixed points $\mp i\sqrt{\alpha}$ is \emph{locally asymptotically stable} respectively.

\subsection{Permanence and Boundedness of Solutions}

\begin{theorem}
  Let $\abs{\alpha}<\abs{\frac{(\gamma-\beta)^2}{4}}$ and the initial values $z_{-k}, z_{-k+1},\dots,z_{-1},z_0 \in B(0,\abs{\bar{z}_1})$. If $\{z_n\}$ is any solution of Eq.(\ref{equation:total-equationA}), then $z_n \in B(0,\abs{\bar{z}_1})$ for all $n\geq 1$, that is $\{z_n\}$ is permanent and bounded.\\
\end{theorem}

\noindent
Proof: If the the initial values $z_{-k}, z_{-k+1},\dots,z_{-1},z_0 \in B(0,\abs{\bar{z}_1})$, then $\abs{z_{-k}}<\abs{\bar{z}_1}, \abs{z_{-k+1}}<\abs{\bar{z}_1},\dots,\abs{z_{-1}}<\abs{\bar{z}_1}$ and $\abs{z_{0}}<\abs{\bar{z}_1}$, we have $$0<\abs{z_1}=\abs{\frac{\alpha+\beta z_{-k}}{\gamma-z_0}}<\abs{\frac{\alpha+\beta \bar{z}_1}{\gamma-\bar{z}_1}}=\abs{\bar{z}_1}$$
Similarly,
$$0<\abs{z_2}=\abs{\frac{\alpha+\beta z_{-k+1}}{\gamma-z_1}}<\abs{\frac{\alpha+\beta \bar{z}_1}{\gamma-\bar{z}_1}}=\abs{\bar{z}_1}$$

\noindent
Using \emph{Mathematical Induction} we have $0<\abs{z_n}<\abs{\bar{z}_1}$ for all $n \geq 1$. Therefore, the open ball $B(0,\abs{\bar{z}_1})$ is invariant of Eq.(\ref{equation:total-equationA}) and the solution $\{z_n\}$ of Eq.(\ref{equation:total-equationA}) is permanent and bounded.

\begin{theorem}
  Let $\abs{\alpha}=\abs{\frac{(\gamma-\beta)^2}{4}}$ and the initial values $z_{-k}, z_{-k+1},\dots,z_{-1},z_0 \in B(0,\abs{\frac{\gamma-\beta}{2}})$ where $\abs{\beta} < \abs{3\gamma}$. If $\{z_n\}$ is any solution of Eq.(\ref{equation:total-equationA}), then $z_n \in B(0,\abs{\frac{\gamma-\beta}{2}})$ for all $n\geq 1$, that is $\{z_n\}$ is \emph{permanent and bounded}.\\
\end{theorem}

\noindent
Proof: Considering $\abs{\beta} < \abs{3\gamma}$, we have $0<\abs{\frac{\gamma-\beta}{2}}<\abs{\gamma}$ and the initial values $z_{-k}, z_{-k+1},\dots,z_{-1},z_0 \in B(0,\abs{\frac{\gamma-\beta}{2}})$, then $\abs{z_{-k}}<\abs{\frac{\gamma-\beta}{2}}, \abs{z_{-k+1}}<\abs{\frac{\gamma-\beta}{2}},\dots,\abs{z_{-1}}<\abs{\frac{\gamma-\beta}{2}}$ and $\abs{z_{0}}<\abs{\frac{\gamma-\beta}{2}}$, we have $$0<\abs{z_1}=\abs{\frac{\alpha+\beta z_{-k}}{\gamma-z_0}}<\abs{\frac{\alpha+\beta \frac{\gamma-\beta}{2}}{\gamma-\frac{\gamma-\beta}{2}}}=\abs{\frac{\gamma-\beta}{2}}$$
Similarly,
$$0<\abs{z_2}=\abs{\frac{\alpha+\beta z_{-k+1}}{\gamma-z_1}}<\abs{\frac{\alpha+\beta \frac{\gamma-\beta}{2}}{\gamma-\frac{\gamma-\beta}{2}}}=\abs{\frac{\gamma-\beta}{2}}$$

\noindent
Using \emph{Mathematical Induction} we have $0<\abs{z_n}<\abs{\frac{\gamma-\beta}{2}}$ for all $n \geq 1$. Therefore, the open ball $B(0,\abs{\frac{\gamma-\beta}{2}})$ is invariant of Eq.(\ref{equation:total-equationA}) and the solution $\{z_n\}$ of Eq.(\ref{equation:total-equationA}) is \emph{permanent and bounded}.\\

\noindent
An example is taken here as an evidence of the \textbf{Theorem 2.4}. Consider $\gamma=5+6i$, $\beta=2+3i$ and so $\alpha=\abs{\frac{(\gamma-\beta)^2}{4}}=0+4.5i$. Here the fixed point is $\frac{\gamma-\beta}{2}=1.5+1.5i$ with modulus  $2.1213$. Also we assume k=2 and the initial values $z_{-2}=i, z_{-1}=1+i$ and $z_0=1.5+i$ note that $\abs{z_{-2}},\abs{z_{-1}}$ and $\abs{z_0}$ are all less than $\abs{\frac{\gamma-\beta}{2}}=2.1213$ and $\abs{\beta}<3\abs{\gamma}$. So the \textbf{Theorem 2.4} is applicable here. What we expect is that all the solution $\{z_n\}$ of Eq.(\ref{equation:total-equationA}) is \emph{permanent and bounded}, that is $0<\abs{z_n}<\abs{\frac{\gamma-\beta}{2}}$ for all $n \geq 1$. The complex trajectory plot including its time series plot are given in the Fig. 5.

\begin{figure}[H]
      \centering

      \resizebox{10cm}{!}
      {
      \begin{tabular}{c c}
      \includegraphics [scale=5]{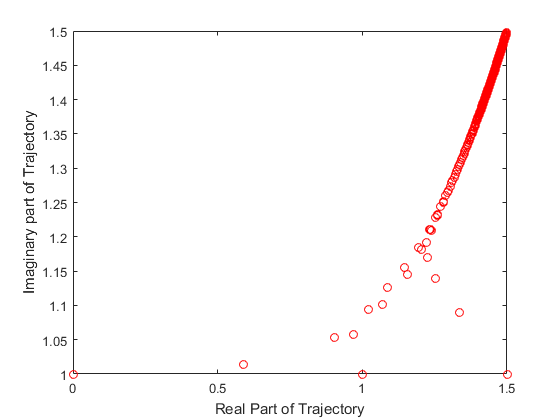}
      \includegraphics [scale=5]{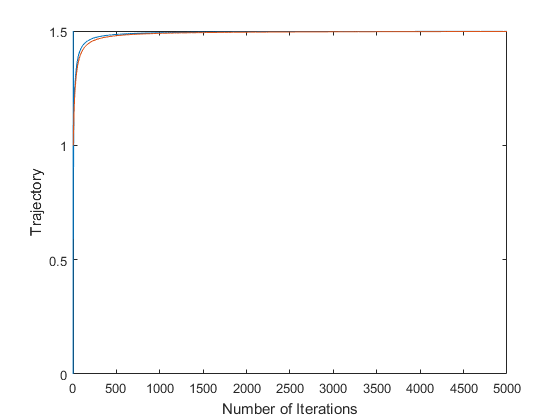}\\
\end{tabular}
      }
\caption{Complex trajectory plot including its time series plot.}

      \end{figure}

 \noindent
 In the Fig. $5$, it is seen that all the solutions $\{z_n\}$ of Eq.(\ref{equation:total-equationA}) are lying in the disk $B(0,\abs{\frac{\gamma-\beta}{2}})$.

\section{Periodic of Solutions}

\noindent
We shall first look for the \emph{prime period two solutions} of the three difference equation Eq.(\ref{equation:total-equationA}).\\

\noindent
Let $\ldots, \phi ,~\psi , ~\phi , ~\psi ,\ldots$, $\phi  \neq \psi $ be a \emph{prime period two solution} of the difference equation Eq.(\ref{equation:total-equationA}). \\
\noindent
If $k$ is even, then $z_n=z_{n-k}$ and then $\psi=\frac{\alpha+\beta\phi}{\gamma-\phi}$ and $\phi=\frac{\alpha+\beta\psi}{\gamma-\psi}$. By solving these two equations we get $(\gamma+\beta) (\psi-\phi) =0$. So $\gamma=-\beta$ will lead to \emph{prime period two solutions} of the difference equation Eq.(\ref{equation:total-equationA}).\\
\noindent
if $k$ is odd, then $z_{n+1}=z_{n-k}$ and then $\psi=\frac{\alpha+\beta\psi}{\gamma-\phi}$ and $\phi=\frac{\alpha+\beta\phi}{\gamma-\psi}$. By solving these two equations we get $(\gamma-\beta) (\psi-\phi) =0$. Therefore $\gamma=\beta$ will lead to \emph{prime period two solutions} of the difference equation Eq.(\ref{equation:total-equationA}).\\

\noindent
It is noted that for $\gamma>\beta>0$ in the real line, there was no \emph{prime period two solution} of the difference equation Eq.(\ref{equation:total-equationA}). Here we list a few examples where we found the higher order periods viz. $9, 10, 13$ and etc.

\begin{table}[H]

\begin{tabular}{| m{5cm} || m{5cm} || m{5cm}|}
\hline
\centering   \textbf{Parameters}: $\alpha$, $\beta$, $\gamma$ \textbf{and Delay term}: $k$ &
\begin{center}
\textbf{Periodic Solutions}
\end{center}
 &
\begin{center}
\textbf{Remarks}
\end{center}\\
\hline
\hline
\centering $\alpha=0.2729 + 0.0372i$, $\beta=0.7690 + 0.3960i$, $\gamma=-\beta$ and $k=1$ &
\begin{center}
$\psi=-0.0444 + 0.1575i$ $\phi=0.2333 + 1.6669i$, Period: 2
\end{center} &
\begin{center}
\includegraphics[width=0.3\textwidth, height=20mm]{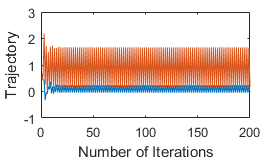}
\end{center}
\\
\hline
\centering $\alpha= 0.4024 + 0.9828i$, $\beta= 0.0740 + 0.6841i$, $\gamma=\beta$ and $k=2$ &
\begin{center}
$\psi=-0.5589 - 0.1317i$ $\phi=-0.1991 + 0.9915i$, Period: 2
\end{center} &
\begin{center}
\includegraphics[width=0.3\textwidth, height=20mm]{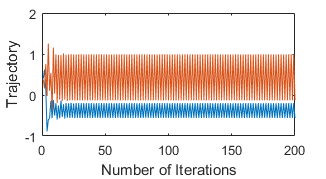}
\end{center}
\\
\hline
\centering $\alpha=-0.1111 + 0.2337i$, $\beta=\alpha$, $\gamma=\beta$ and $k=6$ &
\begin{center}
Period: 9
\end{center}
&
\begin{center}
\includegraphics[width=0.3\textwidth, height=20mm]{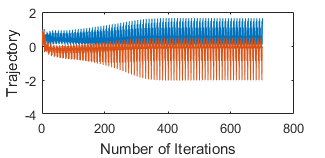}
\end{center}
\\
\hline
\centering $\alpha=-0.1111 + 0.2337i$, $\beta=\alpha$, $\gamma=\beta$ and $k=8$ &
\begin{center}
Period: 11
\end{center}
&
\begin{center}
\includegraphics[width=0.3\textwidth, height=20mm]{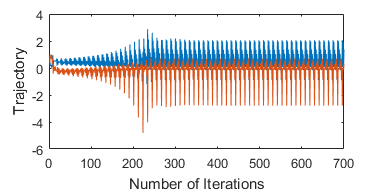}
\end{center}
\\
\hline
\centering $\alpha=-0.1111 + 0.2337i$, $\beta=\alpha$, $\gamma=\beta$ and $k=12$ &
\begin{center}
Period: 15
\end{center}
&
\begin{center}
\includegraphics[width=0.3\textwidth, height=20mm]{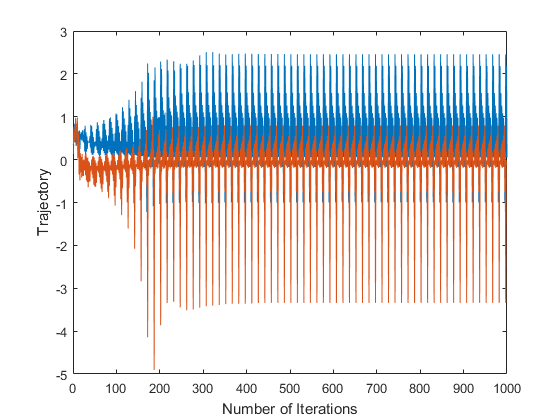}
\end{center}
\\
\hline

\centering $\alpha=-0.1111 + 0.2337i$, $\beta=\alpha$, $\gamma=\beta$ and $k=17$ &
\begin{center}
Period: 20
\end{center}
&
\begin{center}
\includegraphics[width=0.3\textwidth, height=20mm]{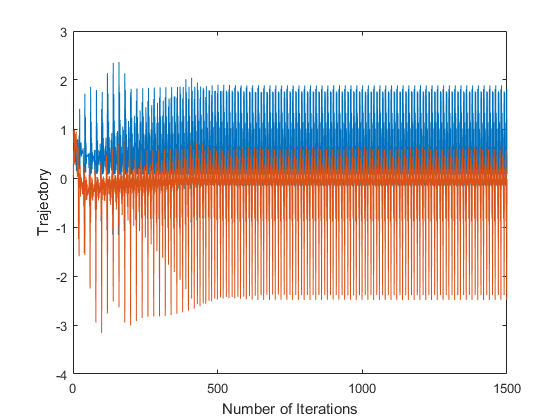}
\end{center}
\\
\hline

\end{tabular}
\caption{Higher Order Periodic Solutions where $\gamma=\pm \beta$, of the equation Eq.(\ref{equation:total-equationA}) for different initial values. In the right most column the corresponding periodic trajectory plots are given.}
\label{Table:}
\end{table}

\begin{table}[H]

\begin{tabular}{| m{5cm} || m{5cm} || m{5cm}|}
\hline
\centering   \textbf{Parameters}: $\alpha$, $\beta$, $\gamma$ \textbf{and Delay term}: $k$ &
\begin{center}
\textbf{Periodic Solutions}
\end{center}
 &
\begin{center}
\textbf{Remarks}
\end{center}\\
\hline
\hline
\centering $\alpha=-0.1111 + 0.2337i$, $\beta=\alpha$, $\gamma=\beta$ and $k=38$ &
\begin{center}
Period: 41
\end{center}
&
\begin{center}
\includegraphics[width=0.3\textwidth, height=20mm]{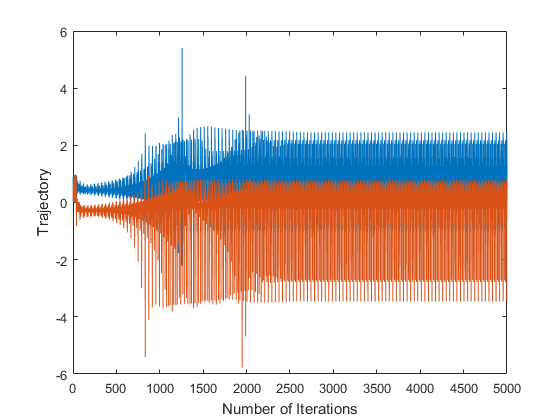}
\end{center}
\\
\hline

\centering $\alpha=-0.1111 + 0.2337i$, $\beta=\alpha$, $\gamma=\beta$ and $k=40$ &
\begin{center}
Period: 43
\end{center}
&
\begin{center}
\includegraphics[width=0.3\textwidth, height=20mm]{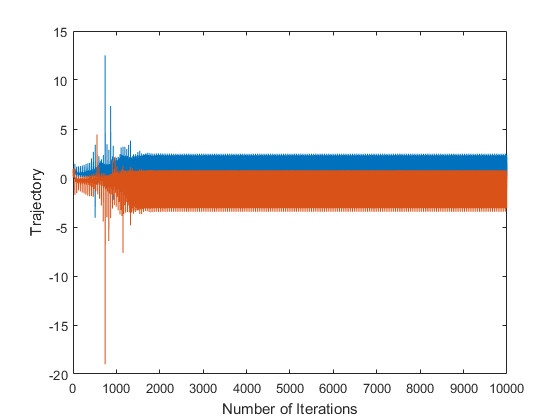}
\end{center}
\\
\hline

\centering $\alpha=-0.1111 + 0.2337i$, $\beta=\alpha$, $\gamma=\beta$ and $k=52$ &
\begin{center}
Period: 55
\end{center}
&
\begin{center}
\includegraphics[width=0.3\textwidth, height=20mm]{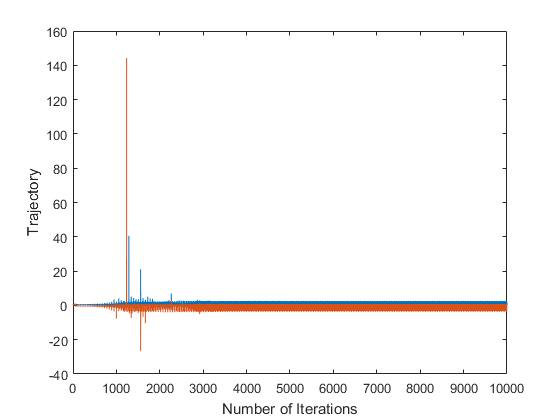}
\end{center}
\\
\hline
\centering $\alpha=-0.1111 + 0.2337i$, $\beta=\alpha$, $\gamma=\beta$ and $k=80$ &
\begin{center}
Period: 83
\end{center}
&
\begin{center}
\includegraphics[width=0.3\textwidth, height=20mm]{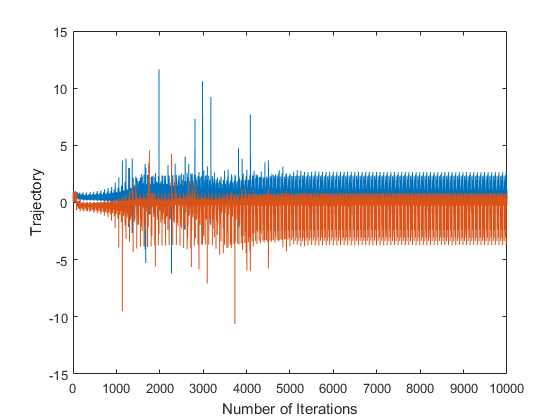}
\end{center}
\\
\hline
\centering $\alpha=-0.1111 + 0.2337i$, $\beta=\alpha$, $\gamma=\beta$ and $k=100$ &
\begin{center}
Period: 103
\end{center}
&
\begin{center}
\includegraphics[width=0.3\textwidth, height=20mm]{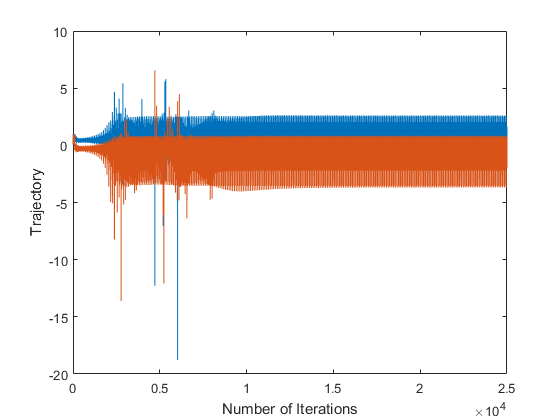}
\end{center}
\\
\hline

\centering $\alpha=-0.1111 + 0.2337i$, $\beta=\alpha$, $\gamma=\beta$ and $k=5100$ &
\begin{center}
Period: 5103
\end{center}
&
\begin{center}
\includegraphics[width=0.3\textwidth, height=20mm]{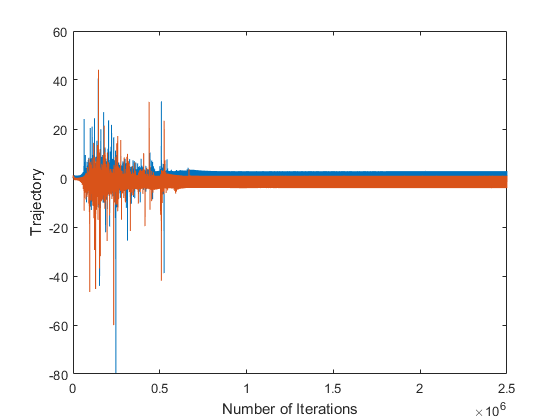}
\end{center}
\\
\hline

\end{tabular}
\caption{Higher Order Periodic Solutions where $\gamma=\pm \beta$, of the equation Eq.(\ref{equation:total-equationA}) for different initial values. In the right most column the corresponding periodic trajectory plots are given.}
\label{Table:}
\end{table}

\noindent
In the Table $2$ and Table $3$, a set of examples of higher order periods are given for different values of $k$ and parameters. What is interesting to note here is that the period is associated to the delay term $k$ as $p=k+3$ ($p$ denotes period) where all the parameters are taken as same ($\alpha=\beta=\gamma$). This inspired us to propose the following conjectures.

\begin{conjecture}
  There exist as many as higher order periods of the difference Eq.(\ref{equation:total-equationA}) is demanded.
\end{conjecture}

\begin{conjecture}
  For $\alpha=\beta=\gamma$, period ($p$) is increased by $3$ with the delay term $k$. i.e. $p=k+3$.
\end{conjecture}

\section{Chaotic Solutions}
A new dynamical behavior of the rational difference equation Eq.(\ref{equation:total-equationA}) is chaoticity which was not present in the real set of parameters. It is really hard to determine the set of all parameters $\alpha$, $\beta$ and $\gamma$ for which the solutions of the Eq.(\ref{equation:total-equationA}) are chaotic but computationally we have encountered some chaotic solutions for some values of the parameters which are given in the following Table. 4. \\
\noindent
   The Lyapunov characteristic exponents serves as a useful tool to quantify chaos. Specifically Lyapunav exponents measure the rates of convergence or divergence of nearby trajectories. Negative Lyapunov exponents indicate convergence, while positive Lyapunov exponents demonstrate divergence and chaos. The magnitude of the Lyapunov exponent is an indicator of the time scale on which chaotic behavior can be predicted or transients decay for the positive and negative exponent cases respectively. In this present study, the largest Lyapunov exponent is calculated for a given solution of finite length numerically \cite{Wolf}.\\
\noindent
From computational evidence, it is arguable that for complex parameters $\alpha$, $\beta$ and $\gamma$ which are stated in the following table the solutions are chaotic for every initial values.

\begin{table}[H]

\begin{tabular}{| m{8cm} | m{6cm} |}
\hline
\centering   \textbf{Parameters} $\alpha$, $\beta$, $\gamma$ and Delay term $k$ &
\begin{center}
\textbf{Lyapunav exponent}
\end{center}\\
\hline
\centering $\alpha=(0.8003, 0.1419)$, $\beta=(0.1576+0.9706)$, $\gamma=0.9572+0.4854)$ and $k=1$ &
\begin{center}
$0.0543$
\end{center}\\
\hline
\centering $\alpha=(0.2217, 0.1174)$, $\beta=(0.6028, 0.7112)$, $\gamma=(0.2217, 0.1174)$ and $k=2$ &
\begin{center}
$1.062$
\end{center}\\
\hline
\centering $\alpha=(0.2564, 0.6135)$, $\beta=(0.6620, 0.4162)$, $\gamma=(0.8419, 0.8329)$ and $k=3$ &
\begin{center}
$1.314$
\end{center} \\
\hline
\centering $\alpha=(0.5349, 0.7210)$, $\beta=(0.4795, 0.6393)$, $\gamma=(0.557, 0.6473)$ and $k=4$ &
\begin{center}
$1.373$
\end{center}\\
\hline

\end{tabular}
\caption{Chaotic solutions of the equation Eq.(\ref{equation:total-equationA} for different choice of parameters and initial values.}
\label{Table:}
\end{table}

\noindent
The chaotic trajectory plots including corresponding complex plots are given the following Fig. 6.

\begin{figure}[H]
      \centering

      \resizebox{11cm}{!}
      {
      \begin{tabular}{c}
      \includegraphics [scale=8]{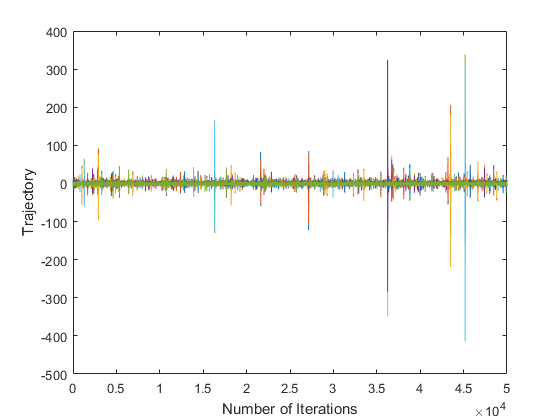}
      \includegraphics [scale=8]{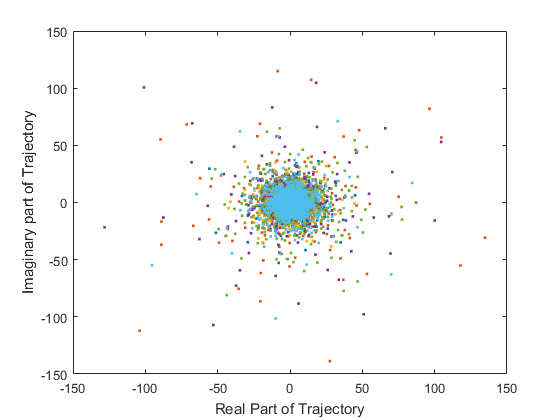}\\
      \includegraphics [scale=8]{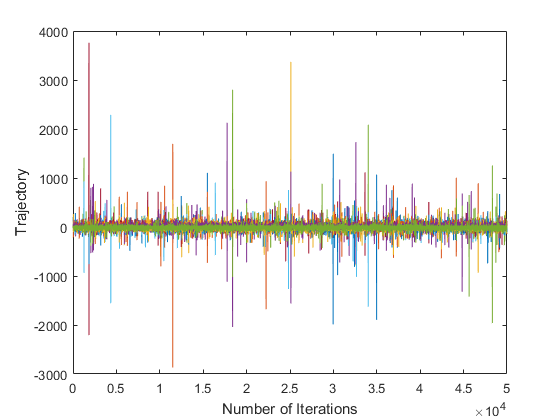}
      \includegraphics [scale=8]{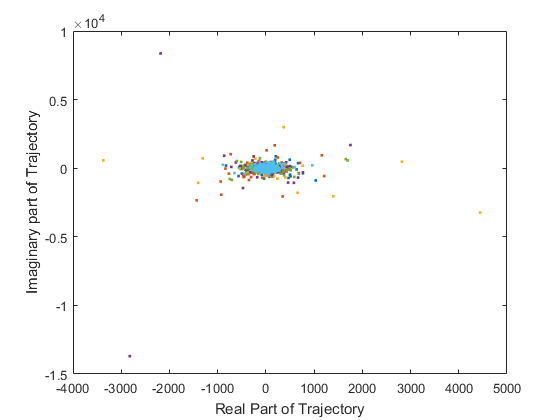}\\
      \includegraphics [scale=8]{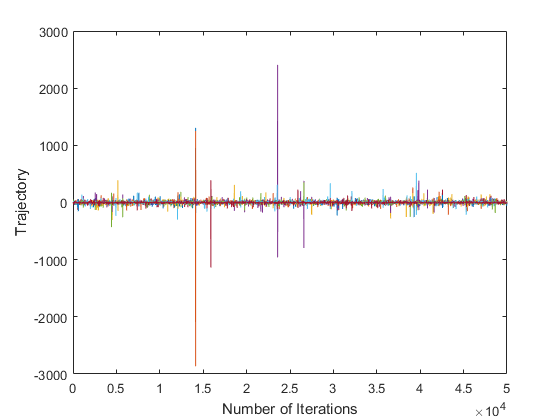}
      \includegraphics [scale=8]{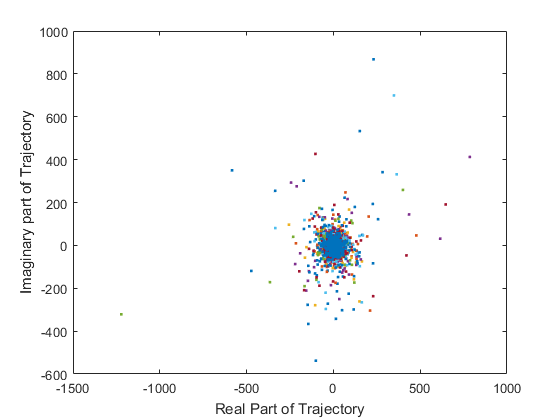}\\
      \includegraphics [scale=8]{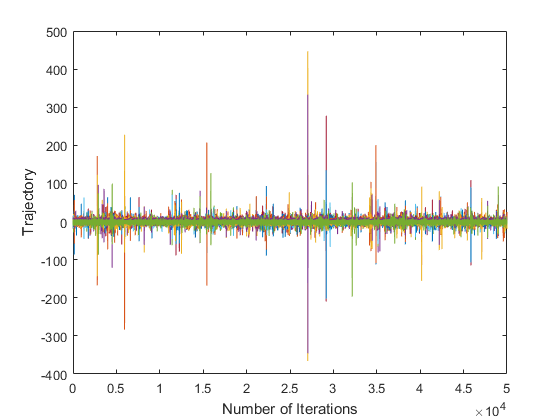}
      \includegraphics [scale=8]{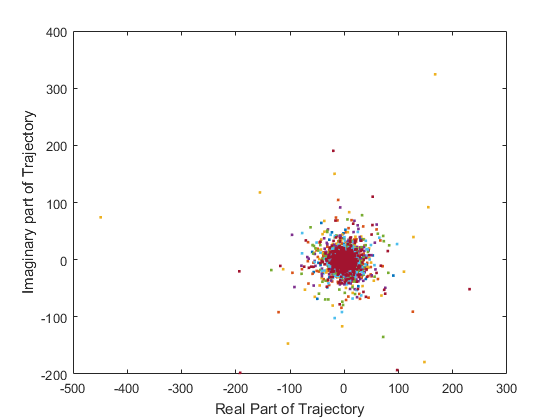}\\
      \end{tabular}
      }
\caption{Chaotic Trajectories of the equation Eq.(\ref{equation:total-equationA} of four different cases as stated in Table 4.}
      \begin{center}

      \end{center}
      \end{figure}

\noindent
In the Fig. 6, for each of the four cases ten different initial values are taken and plotted in the left and in the right corresponding complex plots are given. From the Fig. 6, it is evident that for the four different cases the basin of the chaotic attractor is neighbourhood of the centre $(0, 0)$ of complex plane.
\noindent
It is noted that the \emph{Conjecture} $3.1$ and $3.2$ together suggest that the chaos is simply proportional to the delay term $k$ of the rational dynamical system since the existence of all possible periods suggest that the trajectory is eventually chaotic.

%
%

\end{document}